\documentclass[11pt,openbib]{article}
\usepackage{latexsym}
\usepackage{amsmath}
\usepackage{amssymb}
\usepackage{amsfonts}
\numberwithin{equation}{section}
\newtheorem{thrm}{Theorem}[section]

\newtheorem{defn}{Definition}[section]
\newtheorem{lemm}{Lemma}[section]

\newtheorem{rmrk}{Remark}[section]

\newenvironment{prf}{{\em Proof:}}{\hfill{$\Box$}}
\newcommand{\A}{{\mathcal A}}
\newcommand{\alf}{\alpha}

\newcommand{\dt}[1]{\frac{d}{dt}{#1}}
\newcommand{\dlt}{\delta}
\newcommand{\Dlt}{\Delta}

\newcommand{\dert}{\frac{d}{dt}}

\newcommand{\field}[1]{\mathbb{#1}}
\newcommand{\gmm}{\gamma}
\newcommand{\Gmm}{\Gamma}
\newcommand{\goto}{\rightarrow}
\newcommand{\gs}{\geqslant}

\newcommand{\half}{\frac{1}{2}}
\newcommand{\hf}{\half}

\newcommand{\lmd}{\lambda}

\newcommand{\lf}{\left}
\newcommand{\lra}[2]{\left\langle #1,\,#2\right\rangle}
\newcommand{\ls}{\leqslant}

\newcommand{\mpt}{\mapsto}
\newcommand{\N}{\field{N}}
\newcommand{\nbl}{\nabla}

\newcommand{\ol}{\overline}

\newcommand{\Omg}{\Omega}
\newcommand{\prt}{\partial}
\newcommand{\pder}[2]{\frac{\partial#1}{\partial#2}}
\newcommand{\pdt}[1]{\frac{\partial#1}{\partial t}}
\newcommand{\pdz}[1]{\frac{\partial#1}{\partial z}}

\newcommand{\Di}{\prt_i}

\newcommand{\Px}{\prt_x}
\newcommand{\Py}{\prt_y}
\newcommand{\Pz}{\prt_z}

\newcommand{\Pt}{\partial_t}

\newcommand{\R}{\field{R}}

\newcommand{\ra}{\rightarrow}

\newcommand{\rt}{\right}

\newcommand{\sgm}{\sigma}

\newcommand{\smgrp}{$\{S(t)\}_{t\ge 0}$}

\newcommand{\tht}{\theta}
\newcommand{\ups}{\upsilon}

\newcommand{\veps}{\varepsilon}

\newcommand{\vphi}{\varphi}
\newcommand{\vfi}{\varphi}

\begin{document}
\pagestyle{myheadings}

\title{\bf Finite Dimensionality of the Global Attractor for the Solutions to 3D Primitive Equations with Viscosity}

\author{ Ning Ju $^1$}
\footnotetext[1]{
	Department of Mathematics, Oklahoma State University,
	401 Mathematical Sciences, Stillwater OK 74078, USA.
	Email: {\tt nju@okstate.edu},
	}

\date{July 21, 2015$^2$}
\footnotetext[2]{The original version of this paper was finished on December,
 3rd, 2014. This updated version just includes a short additional
 Section~\ref{s:pbc}.}
\maketitle

\begin{abstract}
 A new method is presented to prove finiteness of the fractal and Hausdorff
 dimensions of the global attractor for the strong solutions to the 3D
 Primitive Equations with viscosity, which is applicable to even more general
 situations than the recent result of \cite{j;t:2014} in the sense that it
 removes {\em all} extra technical conditions imposed by previous analyses.
 More specifically, for finiteness of the dimensions of the global attractor,
 we only need the heat source $Q\in L^2$ which is exactly the condition for
 the existence of global strong solutions and the existence of the global
 attractor of these solutions; while the best existing result, which was
 obtained very recently in \cite{j;t:2014}, still needs the extra condition
 that $\prt_zQ\in L^2$ for finiteness of the dimensions of the global attractor.
 Moreover, the new method can be applied to cases with more complicated
 boundary conditions which present essential difficulties for previous methods.
\\
 


\end{abstract}


\markboth{3D Viscous Primitive Equations}
{\hspace{.2in} N. Ju\hspace{.5in} 3D Viscous Primitive Equations}

\indent
\baselineskip 0.58cm

\section{Introduction}
\label{s:intro}

 Given a bounded domain $D \subset \R^2$ with smooth boundary $\prt D$, we
 consider the following system of viscous Primitive Equations (PEs) of
 Geophysical Fluid Dynamics in the cylinder $\Omg =D\times(-h, 0)\subset\R^3$,
 where $h$ is a positive constant, see e.g. \cite{temma;ziane:04} and the
 references therein:

\noindent
 {\em Conservation of horizontal momentum:}
\begin{equation*}
  \label{e:v}
  \pdt{v} + (v\cdot\nbl)v + w \pdz{v} +\nbl p + fv^\bot + L_1 v = 0;
\end{equation*}
 {\em Hydrostatic balance:}
\begin{equation*}
  \label{e:hs}
  \Pz p + \tht =0 ;
\end{equation*}
 {\em Continuity equation:}
\begin{equation*}
  \label{e:cnt}
  \nbl\cdot v + \Pz w = 0 ;
\end{equation*}
 {\em Heat conduction:}
\begin{equation*}
  \label{e:t}
  \pdt{\tht} + v\cdot\nbl \tht + w \pdz{\tht} + L_2 \tht=Q.
\end{equation*}

 The unknowns in the above system of 3D viscous PEs are the fluid velocity field
 $(v, w)=(v_1, v_2, w)\in\R^3 $ with $v = (v_1, v_2)$ and $v^\bot=(-v_2, v_1)$
 being horizontal, the temperature $\tht$ and the pressure $p$. The Coriolis
 rotation frequency $f = f_0(\beta + y)$ in the $\beta$-plane approximation and
 the heat source $Q$ are given. For the issue concerned in this article, $Q$ is
 assumed to be independent of $t$. In the above equations and in this article,
 $\nbl$ and $\Dlt$ denote the horizontal gradient and Laplacian:
\[ \nbl := (\Px, \Py)  \equiv (\prt_1,\prt_2), \quad
   \Dlt := \Px^2+\Py^2 \equiv \sum_{i=1}^2\Di^2. \]
 The viscosity and the heat diffusion operators $L_1$ and $L_2$ are given
 respectively as follows:
\begin{equation*}
 \label{e:L1.L2}
   L_i := - \nu_i\Dlt - \mu_i \frac{\prt^2}{\prt z^2}, \quad
	i =1, 2,
\end{equation*}
 where the positive constants $\nu_1, \mu_1$ are the horizontal and vertical
 viscosity coefficients and the positive constants $\nu_2, \mu_2$ are the
 horizontal and vertical heat diffusivity coefficients.

 The boundary of $\Omg$ is partitioned into three parts:
 $\partial\Omg = \Gmm_u \cup \Gmm_b \cup \Gmm_s$, where
\begin{align*}
\Gmm_u &:= \{(x, y, z) \in \overline{\Omg} : z = 0\},\\
\Gmm_b &:= \{(x, y, z) \in \overline{\Omg} : z = -h\},\\
\Gmm_s &:= \{(x, y, z) \in \overline{\Omg}:
 (x,y)\in \prt D\}.
\end{align*}

 Consider the following boundary conditions of the PEs as in \cite{cao;titi:05}
 and \cite{j:pe}:
\begin{equation*}
\begin{aligned}
  \mbox{on} \quad \Gmm_u \mbox{\,:}& \quad 
  \pdz{v}= h \tau, \quad w=0, \quad \pdz{\tht}=-\alf(\tht-\Theta),\\
  \mbox{on} \quad \Gmm_b \mbox{\,:}& \quad
  \pdz{v}=0,  \quad w=0, \quad \pdz{\tht}=0,\\
  \mbox{on} \quad \Gmm_s \mbox{\,:}& \quad 
  v\cdot n=0,  \quad \pder{v}{n}\times n= 0,
  \quad \pder{\tht}{n}=0,
\end{aligned}
\end{equation*}
 where $\tau(x, y)$ and $\Theta(x, y)$ are respectively the wind stress and
 typical temperature distribution on the surface of the ocean, $n$ is the
 normal vector of $\Gmm_s$ and $\alf$ is a non-negative constant.
 The above system of PEs will be solved with suitable initial conditions.

 We assume that $Q$, $\tau$ and $\Theta$ are independent of time. Notice that
 results similar to those to be presented here for the autonomous case can
 still be obtained for the non-autonomous case with proper modifications.
 For the autonomous case, assuming some natural compatibility conditions on
 $\tau$ and $\Theta$, one can further set $\tau=0$ and $\Theta=0$ without
 losing generality. See \cite{cao;titi:05} for a detailed discussion on this
 issue.

 Setting $\tau=0$, $\Theta=0$ and using the fact that
\begin{equation*}
 w(x, y, z, t) =  - \int_{-h}^z \nbl\cdot v(x, y, \xi, t)d\xi,
\end{equation*}
\begin{equation*}
 p(x, y, z, t) =  p_s(x, y, t) -\int^z_{-h} \tht(x, y, \xi, t)d\xi,
\end{equation*}
 one obtains the following equivalent formulation of the system of PEs:
\begin{equation}
\label{e:v.n}
\begin{split}
 \pdt{v} + L_1v 
 &+ (v\cdot\nbl)v
 -\left(\int_{-h}^z\nbl\cdot v(x,y,\xi,t)d\xi\right)\pdz{v}\\
 &+ \nbl p_s(x,y,t)-\int_{-h}^z \nbl\tht(x,y,\xi,t)d\xi
 + fv^\bot =0.
\end{split}
\end{equation}
\begin{equation}
\label{e:t.n}
 \pdt{\tht} + L_2\tht + v\cdot\nbl\tht 
 - \left(\int_{-h}^z\nbl\cdot v(x,y,\xi,t)d\xi\right)\pdz{\tht}
 =Q ;
\end{equation}
\begin{equation}
\label{e:bc.v}
 \pdz{v} \Big\vert_{z=0} = \pdz{v}\Big|_{z=-h} =0, \quad
 v\cdot n\big|_{\Gmm_s} =0, \quad \pder{v}{n}\times n\Big|_{\Gmm_s} =0,
\end{equation}
\begin{equation}
\label{e:bc.t}
 \lf(\pdz{\tht}+\alf \tht\rt)\Big|_{z=0}= \pdz{\tht}\Big|_{z=-h}=0, \quad
 \pder{\tht}{n}\Big|_{\Gmm_s} =0,
\end{equation}
\begin{equation}
\label{e:ic.n}
 v(x, y, z, 0) = v_0(x, y, z), \quad \tht(x, y, z, 0) = \tht_0(x, y, z).
\end{equation}

 We remark that the expressions of $w$ and $p$ via integrating the continuity
 equation and the hydrostatic balance equation were already used in
 \cite{lions;temam;wang:92.2} dealing with the Primitive Equations for large
 scale oceans. See also \cite{lions;temam;wang:92.1} for a similar treatment
 of the Primitive Equations for atmosphere.

 Notice that the effect of the {\em salinity} is omitted in the above 3D
 viscous PEs for brevity of presentation. However, our results in this article
 are still valid when the effect of salinity is included. Notice also that
 the right-hand side term of \eqref{e:v.n} is set as $0$. This is just for
 brevity of presentation and it is not technically essential here: if it is
 replaced by a non-zero given external force $R\in L^2(\Omega)$, the results
 of this paper are still valid.

 To the best of our knowledge, the mathematical framework of the viscous
 primitive equations for the large scale ocean was first formulated in
 \cite{lions;temam;wang:92.2}; the notions of weak and strong solutions were
 defined and existence of weak solutions was proved. Uniqueness of weak
 solutions is still unresolved yet.
 Existence of strong solutions {\em local in time} and their uniqueness were
 obtained in \cite{g;m;r:01} and \cite{temma;ziane:04}. Existence of strong
 solutions {\em global in time} was proved independently in \cite{cao;titi:05}
 and \cite{kob:07}. See also \cite{kz:07} for dealing with some other boundary
 conditions. In \cite{j:pe}, existence of the global attractor for the strong
 solutions of the system is proved in the functional space of strong solutions.

 This article focus on the study of finiteness of the dimensions of the global
 attractor for the strong solutions of the system of 3D viscous PEs. This
 result was previously announced in \cite{j:pe}. For simplicity of discussion,
 we set the right-hand side of \eqref{e:v.n} as zero. If this term is replaced
 by a time-independent term $R$, no essential change is needed to be made on
 our analysis. For the main results of this article, one just need to add the
 assumption that $R\in L^2$.

 In a recent work \cite{chu:12}, finiteness of the dimensions of the global
 attractor for the strong solutions for the 3D viscous PEs is obtained under
 {\em periodic} boundary conditions. The proof given in \cite{chu:12} of the
 Ladyzhenskaya squeezing property for the semigroup is based on the fact that,
 for the case of {\em periodic} boundary conditions, the $L^2$ norm of the
 second order spacial derivatives of the solutions is bounded uniformly in
 time and uniformly on the global attractor. This can be proved using previous
 results and analysis for uniform boundedness given in \cite{j:pe} as outlined
 in \cite{chu:12} and \cite{petcu:06}, since one can freely integrate by parts
 without having any boundary term. However, for the case with non-periodic
 boundary conditions such as given by \eqref{e:bc.v} and \eqref{e:bc.t}, some
 complicated boundary terms can not be avoided with the integration by parts
 following the strategy of \cite{chu:12} and \cite{petcu:06}. These boundary
 terms present essential difficulties for {\em a priori} estimates in the $H^2$
 norm. To deal with these difficulties, a more recent work \cite{j;t:2014}
 provides an involved analysis, which proves the existence of a bounded
 absorbing ball in $H^2$ and the uniform boundedness of the $H^2$ norm of the
 solutions. Indeed, before the work of \cite{j;t:2014}, it was not known
 whether or not uniform boundedness for the solutions in $H^2$ is still valid,
 no matter how smooth the initial data and the right-hand side terms are.
 As an application of the uniform $H^2$ boundedness, it is proved in
 \cite{j;t:2014} that the Ladyzhenskaya squeezing property of the semigroup is
 indeed still valid for the solutions on the global attractor in the case of
 the non-periodic boundary conditions as considered here. Thus, the problem of
 finiteness of the dimension of the global attractor for the strong solutions
 of the system \eqref{e:v.n}-\eqref{e:ic.n} measured in the space of strong
 solutions is positively resolved.

 Compared with \cite{chu:12} and \cite{petcu:06}, the analysis of
 \cite{j;t:2014} achieves two improvements. Firstly, the new analysis is
 applicable to both the non-periodic case and the periodic case, while the
 previous methods do not seem to apply to non-periodic boundary conditions,
 such as \eqref{e:bc.v} and \eqref{e:bc.t}. Secondly, the new result of
 \cite{j;t:2014} requires less demanding conditions than in previous works.
 More specifically, \cite{j;t:2014} requires only that $Q, \prt_zQ\in L^2$,
 instead of requiring that $Q\in H^1$ and $\prt_zQ\in L^6$ in addition to
 periodicity as needed in \cite{chu:12}. See Theorem~\ref{t:dim.1} in
 Section~\ref{s:pre} of this article for details.

 However, recalling the main result of \cite{j:pe} (see e.g.
 Theorem~\ref{t:attractor} in Section~\ref{s:pre} of this article), we notice
 that $Q\in L^2$ is enough for the global existence of strong solutions and for
 the existence of the global attractor $\A$ for strong solutions. So, a very
 natural question out of curiosity is that whether or not the additional
 condition $\prt_zQ\in L^2$ imposed in Theorem~\ref{t:dim.1} is essentially
 necessary for the dimensions of the global attractor to be finite. Another
 issue involved with Theorem~\ref{t:dim.1} is that the condition $\alf=0$ seems
 indispensable for the boundary conditions in the analysis of \cite{j;t:2014}.
 Resolving these questions and difficulties is the main concern of this article.

 From the geophysics context and the anisotrophic mathematical structure of
 the 3D Primitives equations, it may seem that the condition $\prt_zQ\in L^2$
 might be quite natural in physics and essential in mathematics. A little
 surprisingly, it will be proved in this article that this condition together
 with the condition $\alf=0$ can indeed be completely dropped out. This main
 result of our current article, Theorem~\ref{t:main}, will be proved in
 Section~\ref{s:dim}. To achieve our goal, a new way is discovered to prove
 finiteness of the Hausdorff and fractal dimensions of the global attractor
 for the strong solutions to the 3D Primitive Equations with viscosity.
 An interesting aspect is that it does {\em not} need the uniform boundedness
 of the $H^2$ norms of the solutions, as required by previous methods. What
 will be used instead are some uniform continuity properties for the solutions
 on the global attrctor $\A$. This new idea helps us to successfully remove the
 extra conditions $\prt_zQ\in L^2$ and $\alf=0$ and thus resolve
 {\em completely} the problem of finiteness of dimensions of the global
 attractor for the strong solutions of the system \eqref{e:v.n}-\eqref{e:ic.n}.
 The new approach presented in this paper has its {\em advantage} over the
 previous one used in \cite{j;t:2014}, as demonstrated further by the addional
 example given in Section~\ref{s:pbc}.

 The rest of this article is organized as follows:

 In Section~\ref{s:pre}, we give the notations, briefly review the background
 results and present the problems to be studied and recall some important facts
 crucial to later analysis. In Section~\ref{s:t}, we prove our first main
 result, Theorem~\ref{t:t}, on the uniform boundedness of $(u_t,\tht_t)$ in
 $L^2$ and the existence of a bounded absorbing ball for $(u_t,\tht_t)$, which
 will be needed to prove our next second main result Theorem~\ref{t:uc}.
 In Section~\ref{s:uc}, we state and prove Theorem~\ref{t:uc} about some uniform
 continuity properties for the solutions on the global attractor, which will be
 crucial for our proof of our final result Theorem~\ref{t:main}.
 In Section~\ref{s:dim}, as an application of the previous results, we finally
 prove Theorem~\ref{t:main}, the main result about the finiteness of the
 Hausdorff and fractal dimensions of the global attractor as measured in $V$,
 the space of strong solutions.
 In Section~\ref{s:pbc}, we briefly mention the case with ``physical boundary
 conditions'' \eqref{e:bc.v.p} on the velocity field $v$, for which the new
 method used in this paper can be easily applied to obtain same conclusions as
 for the case with $v$ satisfying the boundary conditions \eqref{e:bc.v}.
 However, it seems rather difficult for the approach of \cite{j;t:2014} to deal
 with the case of ``physical boundary conditions''.

\section{Preliminaries}
\label{s:pre}

\noindent

 We recall that $D$ is a bounded smooth domain in $\R^2$ and
 $\Omg=D\times[0,-h]$, where $h$ is a positive constant. We denote by
 $L^p(\Omg)$ and $L^p(D)$ ($1\ls p <+\infty$) the classic $L^p$ spaces with
 the norms:
 \begin{equation*}
 \|\phi\|_p = \left\{ \begin{array}{ll}
   \left( \int_\Omg|\phi(x,y,z)|^pdxdydz\right)^\frac{1}{p},
	& \forall \phi\in L^p(\Omg); \\
  \left( \int_D|\phi(x,y)|^pdxdy \right)^\frac{1}{p},
	& \forall \phi\in L^p(D).
  \end{array} \right.
 \end{equation*}
 Denote by $H^m(\Omg)$ and $H^m(D)$ ($m\gs 1$) the classic Sobolev spaces for
 square-integrable functions with square-integrable derivatives up to order $m$.
 We do not distinguish the notations for vector and scalar function spaces,
 which are self-evident from the context. For simplicity, we may use $d\Omg$ to
 denote $dxdydz$ and $dD$ to denote $dxdy$, or we may simply omit them when
 there is no confusion. Using the H\"older inequality, it is easy to show that, 
 for $\vphi\in L^p(\Omg)$,
\begin{equation}
\label{e:lp.av}
 \|\ol{\vphi}\|_{L^p(\Omg)} 
 = 
 h^\frac{1}{p}\|\ol{\vphi}\|_{L^p(D)}
 \ls 
 \|\vphi\|_p, \quad \forall p\in[1, +\infty],
\end{equation}
 where $\ol{\vphi}$ is defined as the vertical average of $\vphi$:
\[ \ol{\vphi}(x,y) = h^{-1}\int_{-h}^0 \vphi(x,y,z)dz .\]

 Define the function spaces $H$ and $V$ as follows:
\begin{equation*}
\begin{split}
 H := H_1\times H_2 := \{ v\in L^2(\Omg)^2\ | \
  \nbl\cdot\ol{v} =0,\quad \ol{v}\cdot n |_{\Gmm_s}=0 \}\times L^2(\Omg),\\
 V := V_1\times V_2 := \{ v\in H^1(\Omg)^2\ | \
  \nbl\cdot\ol{v} =0,\quad v\cdot n |_{\Gmm_s}=0 \} \times H^1(\Omg).
\end{split}
\end{equation*}
 Define the bilinear forms: $a_i: V_i\times V_i \ra \R$, $i=1,2$ as follows:
\begin{equation*}
\begin{aligned}
  a_1(v, u) &= \int_\Omg \left(\nu_1 \nbl v_1\cdot \nbl u_1 
		+\nu_1 \nbl v_2\cdot \nbl u_2
		+ \mu_1 v_z\cdot u_z \right)d\Omg;\\
  a_2(\tht, \eta)& = \int_\Omg \left(\nu_2 \nbl\tht \cdot \nbl\eta 
		+ \mu_2 \tht_z \eta_z \right)d\Omg
                + \alf\int_{\Gmm_u}\tht\eta dxdy.
\end{aligned}
\end{equation*}
 Let $V_i'$ ($i=1,2$) denote the dual space of $V_i$. We define the linear
 operators $A_i : V_i \mapsto V_i'$, $i=1,2$ as follows:
 \begin{equation*}
 \lra{A_1v}{u}=a_1(v, u), \quad \forall v, u \in V_1;
 \quad
 \lra{A_2\tht}{\eta}=a_2(\tht, \eta),
 \quad  \forall \tht, \eta \in V_2,
 \end{equation*}
 where $\lra{\cdot}{\cdot}$ is the corresponding scalar product between $V_i'$
 and $V_i$. We also use $\lra{\cdot}{\cdot}$ to denote the inner products in
 $H_1$ and $H_2$. Define:
 \[ D(A_i) = \{ \phi \in V_i, A_i\phi\in H_i \},\quad i=1,2. \]
 Since $A_i^{-1}$ is a self-adjoint compact operator in $H_i$, by the classic
 spectral theory, the power $A_i^s$ can be defined for any $s\in\R$.
 Then $D(A_i)' = D(A_i^{-1})$ is the dual space of $D(A_i)$ and
 $ V_i = D(A_i^\half)$, $ V_i' = D(A_i^{-\half})$. Moreover,
  \[ D(A_i) \subset V_i \subset H_i \subset V_i' \subset D(A_i)', \]
 where the embeddings above are all compact.
 Define the norm $\|\cdot\|_{V_i}$ by:
 \[ \|\cdot\|_{V_i}^2 = a_i(\cdot,\cdot)
	= \lra{A_i\cdot}{\cdot}
	= \lra{A_i^\half\cdot}{A_i^\half\cdot},
	\quad i=1,2. \]
 The Poincar\'e inequalities are valid. There is a constant $c>0 $,
 such that for any $\phi =(\phi_1, \phi_2)\in V_1$ and $\psi \in V_2$
\begin{equation}
\label{e:pnc}
 c\|\phi\|_2 \ls \|\phi\|_{V_1}, \quad
 c\|\psi\|_2 \ls \|\psi\|_{V_2}.
\end{equation}
 Therefore, there exist constants $c>0$ and $C>0$ such that for any
 $\phi =(\phi_1, \phi_2)\in V_1$ and $\psi \in V_2$,
\begin{equation*}
\label{e:nmeqv}
 c\|\phi\|_{V_1} \ls \|\phi\|_{H^1(\Omg)} \ls C \|\phi\|_{V_1},
\quad
 c\|\psi\|_{V_2} \ls \|\psi\|_{H^1(\Omg)} \ls C \|\psi\|_{V_2}.
\end{equation*}
 Notice that, in the above first inequality, we have written
 $\|\phi\|_{H^1(\Omg)}$ instead of $\|\phi\|_{H^1(\Omg)^2}$. We could also simply
 write $\|\phi\|_{H^1}$. Here and later on as well, we do not distinguish the
 notations for vector and scalar function spaces which are self-evident from
 the context. In this article, we use $c$ and $C$ to denote generic positive
 constants, the values of which may vary from one place to another.

 Recall the following definitions of weak and strong solutions:
\begin{defn}
\label{d:soln}
 Suppose $Q \in L^2(\Omg)$, $(v_0, \tht)\in H$ and $T>0$. The pair
 $(v, \tht)$ is called a {\em weak solution} of the 3D viscous PEs
 (\ref{e:v.n})-(\ref{e:ic.n}) on the time interval $[0, T]$ if it satisfies
 (\ref{e:v.n})-(\ref{e:t.n}) in the weak sense, and also
\begin{equation*}
 (v, \tht) \in C([0, T];H) \cap L^2(0, T; V), \quad
 \prt_t(v,\tht) \in L^1(0, T; V').
\end{equation*}
 Moreover, if $(v_0, \tht_0)\in V$, a weak solution $(v, \tht)$ is called a
 {\em strong solution} of (\ref{e:v.n})-(\ref{e:ic.n}) on the time interval
 $[0, T]$ if, in addition, it satisfies
\begin{equation*}
 (v,\tht) \in C([0, T]; V) \cap L^2(0, T; D(A_1)\times D(A_2)).
\end{equation*}

\end{defn}

 The following theorem on global existence and uniqueness for the strong
 solutions was proved in \cite{cao;titi:05}. See also a related result in
 \cite{kob:07}.

\begin{thrm}
\label{t:strong.g}
 Suppose $Q\in H^1(\Omg)$. Then, for every $(v_0,\tht_0)\in V$ and $T>0$, there
 exists a unique strong solution $(v, \tht)$ on $[0, T]$ to the system of 3D
 viscous PEs, which depends on the initial data continuously in $H$.
\end{thrm}

\begin{rmrk}
 It is easy to see from the proof of Theorem~\ref{t:strong.g} given in
 \cite{cao;titi:05} that the condition $Q\in H^1(\Omg)$ can be relaxed to
 $Q\in L^6(\Omg)$. Notice that there are gaps between Definition \ref{d:soln}
 and Theorem~\ref{t:strong.g} for the condition on $Q$, for the continuity
 of the strong solution with respect to time and for the continuous dependence
 of the strong solution with respect to initial data.
\end{rmrk}

 We now recall the following result proven in \cite{j:pe} for the existence
 of global attractor $\A$ for the strong solutions of the 3D viscous PEs
 (\ref{e:v.n})-(\ref{e:ic.n}).
\begin{thrm}
\label{t:attractor}
 Suppose that $Q\in L^2(\Omg)$ is independent of time. Then the solution
 operator \smgrp \ of the 3D viscous PEs (\ref{e:v.n})-(\ref{e:ic.n}):
 $S(t)(v_0, \tht_0) =(v(t), \tht(t))$ defines a semigroup in the space
 $V$ for $t\in \R_+$.  Moreover, the following statements are valid:
 \begin{enumerate}
 \item
 For any $(v_0, \tht_0)\in V$, $t\mpt S(t)(v_0,\tht_0)$ is continuous from
 $\R_+$ into $V$.
 \item
 For any $t>0$, $S(t)$ is a continuous and compact map in $V$.
 \item
 \smgrp \ possesses a global attractor $\A$ in the space $V$. The global
 attractor $\A$ is compact and connected in $V$ and it is the minimal bounded
 attractor in $V$ in the sense of the set inclusion relation; $\A$ attracts all
 bounded subsets of $V$ in the norm of $V$.
 \end{enumerate}
\end{thrm}

 Recall also the following result proved in \cite{j;t:2014} for finiteness of
 the Hausdorff and fractal dimensions of the global attractor $\A$ as obtained
 in Theorem~\ref{t:attractor}:

\begin{thrm}
\label{t:dim.1}
 Suppose  $\alf=0$ and $Q, Q_z\in L^2(\Omg)$. Then the global attractor
 $\A$ has finite Hausdorff and fractal dimensions measured in the $V$ space.
\end{thrm}

 The main goal of this article is to prove finiteness of the Hausdorff and
 fractal dimensions of the global attractor $\A$ in the space $V$ for any
 $\alf\gs0$ and for any $Q\in L^2(\Omg)$, thus dropping all the {\em extra}
 assumptions imposed in Theorem~\ref{t:dim.1}. The formal statement of this
 main result of this article, Theorem~\ref{t:main}, and its proof will be
 presented in Section~\ref{s:dim}.

 We recall the following lemma which will be of critical usefulness for the
 {\em a priori} estimates in the following sections. See \cite{cao;titi:03},
 and also \cite{j:pe}, for a proof.
\begin{lemm}
\label{l:ju}
 Suppose that $\nbl\ups, \vfi\in H^1(\Omg), \psi\in L^2(\Omg)$.
 Then, there exists a constant $C>0$ independent of $\ups, \vfi, \psi$ and $h$,
 such that
\begin{equation*}
\left|
 \lra{\left(\int_{-h}^z\nbl\cdot \ups(x,y,\xi)d\xi\right) \vfi}{\psi}
 \right|
 \ls C \|\nbl\ups\|_2^{\hf}\left\|\nbl \ups\right\|_{H^1}^{\hf}
	\|\vfi\|_{2}^{\hf}\|\vfi\|_{H^1}^{\hf} \|\psi\|_{2}.
\end{equation*}

\end{lemm}

 Similarly, one can prove the following Lemma~\ref{l:ju.2}, which will be of
 critical usefulness in proving our first main result Theorem~\ref{t:t}.

\begin{lemm}
\label{l:ju.2}
 Suppose that $\nbl\ups, \vfi,\nbl\vfi, \psi,\nbl\psi\in L^2(\Omg)$.
 Then, there exists a constant $C>0$ independent of $\ups, \vfi, \psi$ and $h$,
 such that
\begin{equation*}
\left|
 \lra{\left(\int_{-h}^z\nbl\cdot \ups(x,y,\xi)d\xi\right) \vfi}{\psi}
 \right|
 \ls C \|\nbl\ups\|_2
    \|\vfi\|_2^{\hf}\|\nbl\vfi\|_2^{\hf} \|\psi\|_2^\hf\|\nbl\psi\|_2^\hf.
\end{equation*}

\end{lemm}
 We recall the following formulation of the uniform Gronwall lemma,
 the proof of which can be found in \cite{t.inf}.
\begin{lemm}[Uniform Gronwall Lemma]
\label{l:ug}
   Let $g$, $h$ and $y$ be three non-negative locally integrable
   functions on $(t_0, +\infty)$ such that
   \[ \frac{dy}{dt} \ls gy + h, \qquad \forall t\gs t_0, \]
   and
   \[
   \int_t^{t+r}g(s)ds \ls a_1,
   \qquad 
   \int_t^{t+r}h(s)ds \ls a_2
   \qquad 
   \int_t^{t+r}y(s)ds \ls a_3,
   \qquad
   \forall t\gs t_0,
   \]
   where $r$, $a_1$, $a_2$ and $a_3$ are positive constants.
   Then
   \[
   y(t+r) \ls \left(\frac{a_3}{r}+a_2\right)e^{a_1},
   \qquad
   \forall t\gs t_0.
   \]
\end{lemm}

\section{The bounded absorbing ball for $(\prt_tv,\prt_t\tht)$ in the space
 $H_1\times H_2$}
\label{s:t}

 The existence of a bounded absorbing ball for $(\prt_tv,\prt_t\tht)$ in the
 space $H_1\times H_2$, and the uniform boundedness of $\|\prt_tv\|_2$ and
 $\|\prt_t\tht\|_2$ for $t\in[0,+\infty)$, has been proved in \cite{j;t:2014}
 for the case with $\alf=0$ and $Q, \prt_zQ\in L^2$. In the following, we prove
 the same result for the case $\alf\gs0$ under the condition that $Q\in L_2$
 only. Notice that the extra condition $\prt_zQ\in L^2$ is not needed in the
 following Theorem~\ref{t:t}, which is our first main result. It will be used
 for the proof of our final main result of this article.
\begin{thrm}
\label{t:t}
 Suppose $Q \in L^2(\Omg)$ and $\alf\gs0$.

 For any $(v_0,\tht_0)\in V_1\times V_2$ and
 $(\prt_tv(0),\prt_t\tht(0))\in H_1\times H_2$, there exists a unique solution
 $(v,\tht)$ of \eqref{e:v.n}-\eqref{e:ic.n} such that
\begin{equation*}
 (\prt_tv,\prt_t\tht)\in L^\infty(0, +\infty; H_1\times H_2).
\end{equation*}
 Moreover, there exists a bounded absorbing ball for $(\prt_tv,\prt_t\tht)$ in
 the space of $H_1\times H_2$.

\end{thrm}

\begin{prf}

 We prove the case $\alf>0$. The case $\alf=0$ is similar.

 We only need to prove the existence of a bounded absorbing ball for
 $(\prt_tv,\prt_t\tht)$ in the space of $H_1\times H_2$. The uniform boundedness
 of $\|\prt_tv\|_2,\|\prt_t\tht\|_2$ for $t$ in $[0,+\infty)$ will then follow
 easily and the uniqueness of the solution is obvious.

 Denote
\[ u:=v_t=\prt_t v, \quad \zeta:=\tht_t=\prt_t\tht. \]
 Notice that $u$ and $\zeta$ in the proof of Theorem~\ref{t:t} are different
 from those in the proof of Theorem~\ref{t:main} of Section~\ref{s:dim}.

 The proof of Theorem~\ref{t:t} is divided into two steps.

{\bf Step 1}.
 We prove that the time average of $\|(v_t,\tht_t)\|_2^2$ is uniformly bounded
 with respect to $t$.

 Taking the inner product of \eqref{e:v.n} with $v_t$ and using the boundary
 conditions \eqref{e:bc.v}, we find
\begin{equation}
\label{e:v_t}
\begin{split}
 \hf\dert\|v\|_{V_1}^2 + \|u\|_2^2
 =& -\lra{v\cdot\nbl v}{u} -\lra{wv_z}{u}\\
 &+ \lra{\int_{-h}^z\nbl\tht(x,y,\xi,t)d\xi}{u}- \lra{fv^\bot}{u},
\end{split}
\end{equation}
 where we have used the following calculations
\[ \int_\Omg v_t\cdot\Dlt v
 = \int_{-h}^0\lf[\int_{\prt D} v_t\cdot \frac{\prt v}{\prt n}
  -\hf\int_D\prt_t(|\nbl v|^2)\rt] = -\hf\frac{d}{dt}\|\nbl v\|_2^2, \]
\[ \int_\Omg v_t\cdot\prt_z^2 v
 =\int_D\lf[v_t\cdot v_z\big|_{z=-h}^0-\int_{-h}^0\prt_t(v_z)\cdot v_z\rt]
 = -\hf \frac{d}{dt}\|v_z\|_2^2, \]
and
\begin{equation*}
\begin{split}
 \int_\Omg \nbl p_s(x,y,t) \cdot v_t &= \int_{-h}^0
\lf[ \int_{\prt D} p(x,y,t) v_t\cdot n -\int_D p_s(x,y,t) \nbl\cdot v_t\rt] \\
 &= - \int_D p_s(x,y,t) \lf(\int_{-h}^0\nbl\cdot v(x,y,z,t)dz\rt)_t dxdy\\
 &= \int_D p_s(x,y,t) w_t(x,y,0,t) dxdy =0,
\end{split}
\end{equation*}
 since $w$ satisfies the boundary condition $w(x,y, 0, t) =0$ for $(x,y)\in D$.

 From \eqref{e:v_t}, we easily obtain
\begin{equation*}
\begin{split}
 \hf\dert\|v\|_{V_1}^2 &+ \|u\|_2^2\\
\ls& \|v\|_6\|\nbl v\|_3\|u\|_2 + \|w v_z\|_2\|u\|_2
 + C(\|\nbl\tht\|_2+\|v\|_2)\|u\|_2\\
\ls& C_\veps(\|v\|_{H^1}^2\|\nbl v\|_2\|\nbl v\|_{H^1} + \|wv_z\|_2^2
 + \|\nbl\tht\|_2^2 +\|v\|_2^2) + \veps\|u\|_2^2\\
\ls&  C_\veps(\|v\|_{H^1}^2\|\nbl v\|_2\|\nbl v\|_{H^1}
 + \|\nbl v\|_2\|\nbl v\|_{H^1}\|v_z\|_2\|v_z\|_{H^1})\\
& +C_\veps(\|\nbl\tht\|_2^2+\|v\|_2^2) + \veps\|u\|_2^2,
\end{split}
\end{equation*}
 where in the last inequality above, Lemma~\ref{l:ju} is used. Choosing
 $\veps=\hf$, we obtain
\begin{equation}
\label{e:vt.a}
 \dert\|v\|_{V_1}^2 + \|u\|_2^2 \ls C h_1(t),
\end{equation}
 with
\[ h_1(t):= \|v\|_{H^1}^2\|\nbl v\|_2\|\nbl v\|_{H^1}
 + \|\nbl v\|_2\|\nbl v\|_{H^1}\|v_z\|_2\|v_z\|_{H^1}
 + \|\nbl\tht\|_2^2+\|v\|_2^2.\]
 Taking the inner product of \eqref{e:t.n} with $\tht_t$ and using
 \eqref{e:bc.t} yields
\begin{equation}
\label{e:t_t}
\begin{split}
 \frac{\nu_2}{2}\dert\|\nbl\tht\|_2^2 +\frac{\mu_2}{2}\dert\|\tht_z\|_2^2
 &+\frac{\mu_2\alf}{2}\dert\|\tht(z=0)\|_2^2 + \|\zeta\|_2^2\\
=&-\lra{v\cdot\nbl\tht}{\zeta} -\lra{w\tht_z}{\zeta} +\lra{Q}{\zeta},
\end{split}
\end{equation}
 where we have used the following calculations:
\[ \int_\Omg \tht_t \Dlt\tht
 = \int_{-h}^0\lf[\int_{\prt D} \tht_t \frac{\prt\tht}{\prt n}
  -\hf\int_D\prt_t(|\nbl\tht|^2)\rt] = -\hf\frac{d}{dt}\|\nbl\tht\|_2^2, \]
\[ \int_\Omg \tht_t \prt_z^2\tht
 =\int_D\lf[\tht_t \tht_z\big|_{z=-h}^0-\int_{-h}^0\prt_t(\tht_z) \tht_z\rt]
 = -\alf \dert\|\tht(z=0)\|_2^2 -\hf \frac{d}{dt}\|\tht_z\|_2^2. \]
 By \eqref{e:t_t}, we obtain
\begin{equation*}
\begin{split}
 \hf\dert\|\tht\|_{V_2}^2 &+ \|\zeta\|_2^2\\
\ls & \|v\|_6\|\nbl\tht\|_3\|\zeta\|_2 + \|w \tht_z\|_2\|\zeta\|_2 +\|Q\|_2\|\zeta\|_2\\
\ls & C(\|v\|_6^2\|\nbl\tht\|_3^2 +\|w\tht_z\|_2^2+\|Q\|_2^2) +\hf\|\zeta\|_2^2\\
\ls & C(\|v\|_{H^1}^2\|\nbl\tht\|_2\|\nbl\tht\|_{H^1}
 +\|\nbl v\|_2\|\nbl v\|_{H^1}\|\tht_z\|_2\|\tht_z\|_{H^1}+\|Q\|_2^2)
 +\hf\|\zeta\|_2^2,
\end{split}
\end{equation*}
 where in the last inequality we have used Lemma~\ref{l:ju}. Therefore,
\begin{equation}
\label{e:tt.a}
 \dert\|\tht\|_{V_2}^2 + \|\zeta\|_2^2 \ls C h_2(t),
\end{equation}
 with
\[ h_2(t):= \|v\|_{H^1}^2\|\nbl\tht\|_2\|\nbl\tht\|_{H^1}
 + \|\nbl v\|_2\|\nbl v\|_{H^1}\|\tht_z\|_2\|\tht_z\|_{H^1} +\|Q\|_2^2.\]
 Integrating \eqref{e:vt.a} and \eqref{e:tt.a} with respect to $t$ yields
\begin{equation}
\label{e:vt.b}
 \|v(t+1)\|_{V_1}^2 + \int_t^{t+1}\|u(\tau)\|_2^2d\tau
 \ls \|v(t)\|_{V_1}^2 + C\int_t^{t+1}h_1(\tau)d\tau,
\end{equation}
\begin{equation}
\label{e:tt.b}
 \|\tht(t+1)\|_{V_2}^2 + \int_t^{t+1}\|\zeta(\tau)\|_2^2d\tau
 \ls \|\tht(t)\|_{V_2}^2 + C\int_t^{t+1}h_2(\tau)d\tau.
\end{equation}
 Notice that the previous uniform {\em a priori} estimates in \cite{j:pe} yield
 uniform boundedness of
\[ \|v(t)\|_{V_1}, \quad \|\tht(t)\|_{V_2},\quad
\int_t^{t+1}h_1(\tau)d\tau\ \text{ and }\ \int_t^{t+1}h_2(\tau)d\tau, \]
 with respect to $t>0$ and a bounded absorbing set in $\R_+$ for each of the
 above four terms. With these uniform estimates, we can conclude Step 1 from
 \eqref{e:vt.b} and \eqref{e:tt.b}.
 
{\bf Step 2}.
 We now prove the existence of a bounded absorbing ball for $(v_t,\tht_t)$ in
 $H$ and the uniform boundedness of $\|(v_t,\tht_t)\|_2$ for $t\in[0,+\infty)$.
 
 By \eqref{e:v.n} and \eqref{e:t.n}, we have
\begin{equation}
\label{e:vt}
\begin{split}
 u_t + L_1u &+ (u\cdot\nbl)v +(v\cdot\nbl)u  + w_tv_z + w u_z\\
& + \nbl (p_s)_t -\int_{-h}^z\nbl\zeta(x,y,\xi,t)d\xi + fu^\bot =0,
\end{split}
\end{equation}
\begin{equation}
\label{e:tt}
\zeta_t + L_2\zeta + u\cdot\nbl\tht + v\cdot\nbl\zeta
 +w_t\tht_z + w \zeta_z =0.
\end{equation}
 Taking the inner product of \eqref{e:vt} with $u$ and using the boundary
 conditions \eqref{e:bc.v} and \eqref{e:bc.t}, we obtain
\begin{equation}
\label{e:vt.1}
\begin{split}
\hf\dert\|u\|_2^2 +& \nu_1 \|\nbl u\|_2^2 + \mu_1 \| u_z\|_2^2\\
=& -\lra{(u\cdot\nbl)v}{u} -\lra{(v\cdot\nbl)u}{u}
  - \lra{w_tv_z}{u} - \lra{w u_z}{u} \\
 & - \lra{\nbl (p_s)_t}{u} + \lra{\int_{-h}^z\nbl\zeta}{u} - \lra{fu^\bot}{u}\\
=& -\lra{(u\cdot\nbl)v}{u} - \lra{w_tv_z}{u} + \lra{\int_{-h}^z\nbl\zeta}{u}
=: I_1 + I_2 + I_3,
\end{split}
\end{equation}
 where we have used the fact that
\[ \lra{\nbl (p_s)_t}{u} = \lra{fu^\bot}{u}
   = \lra{(v\cdot\nbl)u}{u} + \lra{w u_z}{u} =0. \]
 Taking the inner product of \eqref{e:tt} with $\zeta$ and using \eqref{e:bc.t},
 we obtain
\begin{equation}
\label{e:tt.1}
\begin{split}
 \hf\dert\|\zeta\|_2^2 +& \nu_2\|\nbl\zeta\|_2^2 +\mu_2\|\zeta_z\|_2^2
 +\mu_2\alf\|\zeta(z=0)\|_2^2\\
=& -\lra{u\cdot\nbl\tht}{\zeta} -\lra{v\cdot\nbl\zeta}{\zeta}
 - \lra{w_t\tht_z}{\zeta} -\lra{w \zeta_z}{\zeta}\\
=& -\lra{u\cdot\nbl\tht}{\zeta} - \lra{w_t\tht_z}{\zeta}
=: J_1 + J_2,
\end{split}
\end{equation}
 where we have used the fact that
\[ \lra{v\cdot\nbl\zeta}{\zeta} + \lra{w \zeta_z}{\zeta} =0. \]
 Notice that first,
\begin{equation*}
\begin{split}
 |I_1| \ls& \|\nbl v\|_2\|u\|_4^2
\ls C\|\nbl v\|_2 \|u\|_2^\hf \|u\|_{H^1}^{\frac{3}{2}}\\
\ls& C\|\nbl v\|_2 \|u\|_2^\hf (\|u\|_2^{\frac{3}{2}}
 +\|\nbl u\|_2^{\frac{3}{2}}+\|u_z\|_2^{\frac{3}{2}})\\
\ls & (C\|\nbl v\|_2 +C_\veps\|\nbl v\|_2^4)\|u\|_2^2
 + \veps(\|\nbl u\|_2^2 + \|u_z\|_2^2),
\end{split}
\end{equation*}
 and similarly,
 \begin{equation*}
\begin{split}
|J_1| \ls& \|\nbl\tht\|_2\|u\|_4\|\zeta\|_4
\ls \hf\|\nbl\tht\|_2(\|u\|_4^2+\|\zeta\|_4^2)\\
 \ls& (C\|\nbl\tht\|_2 +C_\veps\|\nbl\tht\|_2^4)(\|u\|_2^2+\|\zeta\|_2^2)\\
  &+ \veps(\|\nbl u\|_2^2 + \|u_z\|_2^2+\|\nbl\zeta\|_2^2 + \|\zeta_z\|_2^2).
\end{split}
\end{equation*}
 Secondly, by Lemma~\ref{l:ju.2}, we have
\begin{equation*}
\begin{split}
|I_2| = & \lf|\lra{\lf(\int_{-h}^z\nbl\cdot u\rt)v_z}{u}\rt|\\
\ls & C\|\nbl u\|_2 \|v_z\|_2^\hf\|\nbl v_z\|_2^\hf\|u\|_2^\hf\|\nbl u\|_2^\hf\\
 =& C\|v_z\|_2^\hf\|\nbl v_z\|_2^\hf\|u\|_2^\hf\|\nbl u\|_2^{\frac{3}{2}}\\
\ls & C_\veps\|v_z\|_2^2\|\nbl v_z\|_z^2\|u\|_2^2 + \veps\|\nbl u\|_2^2,
\end{split}
\end{equation*}
 and similarly,
 \begin{equation*}
\begin{split}
|J_2| =& \lf| \lra{\lf(\int_{-h}^z\nbl\cdot u\rt)\tht_z}{\zeta} \rt|\\
 \ls& C\|\nbl u\|_2 \|\tht_z\|_2^\hf\|\nbl\tht_z\|_2^\hf\|\zeta\|_2^\hf
       \|\nbl\zeta\|_2^\hf\\
\ls& C_\veps\|\tht_z\|_2^2\|\nbl\tht_z\|_2^2\|\zeta\|_2^2
 + \veps \|\nbl\zeta\|_2^2 + \veps\|\nbl u\|_2^2.
\end{split}
\end{equation*}
 Finally, we find
\begin{equation*}
\begin{split}
|I_3| =& \Big|\lra{\int_{-h}^z\nbl\zeta}{u}\Big|
 =\Big|\lra{\int_{-h}^z\zeta}{\nbl\cdot u}\Big|\\
& \ls \|\zeta\|_2\|\nbl u\|_2
\ls C_\veps\|\zeta\|_2^2 + \veps \|\nbl u\|_2^2.
\end{split}
\end{equation*}
 Now, inserting the above estimates on the $I_i$'s and $J_i$'s into
 \eqref{e:vt.1} and \eqref{e:tt.1} and choosing $\veps>0$ sufficiently small,
 we obtain
\begin{equation}
\label{e:vt.tt}
 y'(t) + \gmm(\|u\|_{V_1}^2+\|\nbl_3\zeta\|_{V_2}^2) \ls C g(t) y(t),
\end{equation}
 where $\gmm=\min\{\nu_1,\nu_2,\mu_1,\mu_2\}$ and
\begin{equation*}
\begin{split}
 y(t) &:= \|u(t)\|_2^2 + \|\zeta(t)\|_2^2\\
 g(t) &:=  1 + \|\nbl v\|_2^4 + \|\nbl\tht\|_2^4 +\|v_z\|_2^2\|\nbl v_z\|_2^2
 + \|\tht_z\|_2^2\|\nbl \tht_z\|_2^2.
\end{split}
\end{equation*}

 Now, we can apply the Uniform Gronwall Lemma to \eqref{e:vt.tt} and using the
 previous {\em a priori} uniform estimates proved in \cite{j:pe} and in Step 1,
 we obtain the existence of a bounded absorbing set for $y$ in $\R_+$ and the
 uniform boundedness of $y(t)$ for $t\in\R_+$. This finishes Step 2.

\end{prf}

\section{Uniform Continuity on the Global Attractor}
\label{s:uc}

 In this section,  as an application of Theorem~\ref{t:t}, we prove the
 following theorem about a few important uniform continuity properties for the
 solutions on the global attractor, which will be very crucial in the analysis
 in Section~\ref{s:dim}. Notice that for simplicity of presentation, we do not
 formulate the theorem in its most complete or sharpest form. We just present
 the version of the results sufficient for our purpose in proving our final
 main result, Theorem~\ref{t:main}.

\begin{thrm}
\label{t:uc}
Suppose $Q\in L^2$ and $\alf\gs0$.

Let $(v_0,\tht_0)\in\A$. Then there is a constant $C>0$, which is independent
of $(v_0,\tht_0)$ and $T>0$, such that
\[ \lf| \|v(T)\|_{V_1}^2 - \|v(0)\|_{V_1}^2 \rt| \ls CT^\hf, \quad
   \lf| \|\tht(T)\|_{V_2}^2 - \|\tht(0)\|_{V_2}^2 \rt| \ls CT^\hf, \]
\[ \int_0^T\lf(\|v_z(t)\|_{V_1}^2+\|\nbl v(t)\|_{V_1}^2\rt)dt \ls C(T^\hf + T), \]
 and
\[ \int_0^T\lf(\|\Dlt\tht(t)\|_2^2+\|\tht_z(t)\|_{V_2}^2\rt)dt \ls C(T^\hf + T). \]

\end{thrm}

\begin{prf}

 It should be pointed out that it is proved in \cite{j:pe} that, for any
 $(v_0,\tht_0)\in V_1\times V_2$, the above two integrals are indeed uniformly
 bounded with respect to $T$, and if $(v_0,\tht_0)\in\A$ the bounds of the
 above two integrals are independent of both $T>0$ and initial data. What we
 are interested to here is that, as $T\goto 0+$, these integrals are not only
 continuous with respect to $T$, i.e. they go to zero as $T$ goes to $0$; they
 are indeed {\em uniformly} continuous with respect to $T$ and {\em initial
 data}.

 Notice first that, by Theorem~\ref{t:attractor} and a well known lemma in
 \cite{t:nse}, we have 
\begin{equation*}
\begin{split}
 \|v(T)\|_{V_1}^2 - \|v(0)\|_{V_1}^2
 =& \|A_1^\hf v(T)\|_2^2-\|A_1^\hf v(0)\|_2^2\\
 =&2 \int_0^T\lra{(A_1^\hf v)_t(t)}{A_1^\hf v(t)} dt\\
 =&2 \int_0^T\lra{v_t(t)}{A_1 v(t)} dt.
\end{split}
\end{equation*}
 Therefore,
\begin{equation*}
\begin{split}
\lf|\|v(T)\|_{V_1}^2 - \|v(0)\|_{V_1}^2\rt|
 \ls& 2 \int_0^T\|v_t(t)\|_2\|A_1 v(t)\|_2 dt\\
 \ls& 2 \sup_{t\in(0,T)}\|v_t(t)\|_2\int_0^T\|A_1 v(t)\|_2 dt\\
 \ls& 2\sup_{t\in(0,T)}\|v_t(t)\|_2 T^\hf
      \lf(\int_0^T\|A_1 v(t)\|_2^2 dt\rt)^\hf\\
 \ls& CT^\hf.
\end{split}
\end{equation*}
 Similarly, we have
\begin{equation*}
\begin{split}
\lf|\|\tht(T)\|_{V_2}^2 - \|\tht(0)\|_{V_2}^2\rt|
 \ls& 2\sup_{t\in(0,T)}\|\tht_t(t)\|_2 T^\hf
      \lf(\int_0^T\|A_2 \tht(t)\|_2^2 dt\rt)^\hf\\
 \ls& CT^\hf.
\end{split}
\end{equation*}
 Next, from \cite{cao;titi:05}, we know that
\begin{equation*}
 \frac{d}{dt}\|v_z\|_2^2 +\nu_1\|\nbl v_z\|^2 + \nu \|v_{zz}\|_2^2
 \ls C(\|v\|_6^4\|v_z\|_2^2 +\|\tht\|_2^2).
\end{equation*}
 Therefore,
\begin{equation*}
\begin{split}
 \int_0^T \|v_z(t)\|_{V_1}^2dt
 \ls& \lf|\|v_z(T)\|_2^2-\|v_z(0)\|_2^2\rt|\\
    &+ C \int_0^T (\|v(t)\|_6^4\|v_z(t)\|_2^2+\|\tht(t)\|_2^2dt\\
\ls & CT^\hf + CT.
\end{split}
 \end{equation*}
 Similarly, from \cite{cao;titi:05}, we know that
\begin{equation*}
 \frac{d}{dt}\|\nbl v\|_2^2 + \|\nbl v\|_{V_1}^2 \ls
 C(\|v\|_6^4\|\nbl v\|_2^2 + \|\nbl v_z\|_2^2\|v_z\|_2^2\|\nbl v\|_2^2
   +\|\nbl\tht\|_2^2)
\end{equation*}
 Note that in \cite{cao;titi:05}, the second term on the right-hand side
 of the above inequality is incorrectly written with $\|\nbl v_z\|_2^2$ being
 replaced by $\|\nbl v\|_2^2$.

 Therefore,
\begin{equation*}
\begin{split}
\int_0^T\|\nbl v(t)\|_{V_1}^2dt 
\ls & \lf|\|\nbl v(T)\|_2^2-\|\nbl v(0)\|_2^2\rt|\\
 &+ C\int_0^T(\|v\|_6^4\|\nbl v\|_2^2 + \|\nbl v_z\|_2^2\|v_z\|_2^2\|\nbl v\|_2^2
   +\|\nbl\tht\|_2^2)dt\\
 \ls & CT^\hf + CT +C\int_0^T\|\nbl v_z\|_2^2dt\\
 \ls & CT^\hf + CT. 
\end{split}
\end{equation*}
 Finally, by \eqref{e:t.n} and \eqref{e:bc.t}, we have
\begin{equation*}
\begin{split}
 \hf\dert&\lf(\|\tht_z\|_2^2+\|\nbl\tht\|_2^2+\alf\|\tht(z=0)\|_2^2\rt)\\
 &+\nu_2\|\Dlt\tht\|_2^2+(\nu_2+\mu_2)\|\nbl\tht_z\|_2^2+\mu_2\|\tht_{zz}\|_2^2\\
 &+\alf(\nu_2+\mu_2)\|\nbl\tht(z=0)\|_2^2\\
 =&
 \lra{v\cdot\nbl\tht}{\Dlt\tht+\tht_{zz}} +
 \lra{w\prt_z\tht}{\Dlt\tht+\tht_{zz}}+
 \lra{Q}{\Dlt\tht+\tht_{zz}}\\
 =:& I_1 + I_2 + I_3.
\end{split}
\end{equation*}
 We estimate the $I_i$'s, $i=1,2,3$, as following.
\begin{equation*}
\begin{split}
 |I_1| &\ls \|v\|_6\|\nbl\tht\|_3(\|\Dlt\tht\|_2+\|\tht_{zz}\|_2)\\
 &\ls C\|v\|_6^2\|\nbl\tht\|_3^2 + \veps(\|\Dlt\tht\|_2^2+\|\tht_{zz}\|_2^2)\\
 &\ls C\|v\|_6^2\|\nbl\tht\|_2\|\nbl\tht\|_{V_2}
 +\veps(\|\Dlt\tht\|_2^2+\|\tht_{zz}\|_2^2).
\end{split}
\end{equation*}
By Lemma~\ref{l:ju},
\begin{equation*}
\begin{split}
 |I_2| \ls& C\|\nbl v\|_2^\hf\|\nbl^2v\|_2^\hf
           \|\prt_z\tht\|_2^\hf\|\nbl\prt_z\tht\|_2^\hf
           (\|\Dlt\tht\|_2+\|\tht_{zz}\|_2)\\
 \ls& C\|\nbl v\|_2\|\nbl^2v\|_2
           \|\prt_z\tht\|_2\|\nbl\prt_z\tht\|_2
           + \veps(\|\Dlt\tht\|_2^2+\|\tht_{zz}\|_2^2)\\
 \ls& C\|\nbl v\|_2^2\|\prt_z\tht\|_2^2\|\nbl^2v\|_2^2
 +\veps(\|\nbl\prt_z\tht\|_2^2+\|\Dlt\tht\|_2^2+\|\tht_{zz}\|_2^2).
\end{split}
\end{equation*}
 Finally,
\begin{equation*}
 |I_3| \ls \|Q\|_2(\|\Dlt\tht\|_2+\|\tht_{zz}\|_2)
 \ls C\|Q\|_2^2 + \veps(\|\Dlt\tht\|_2^2+\|\tht_{zz}\|_2^2).
\end{equation*}
 Hence,
\begin{equation*}
\begin{split}
 \dert\|\tht\|_{V_2}^2 &+ \|\Dlt\tht\|_2^2 + \|\tht_z\|_{V_2}^2\\
 \ls& C(\|v\|_6^2\|\nbl\tht\|_2\|\nbl\tht\|_{V_2}
 +\|\nbl v\|_2^2\|\tht_z\|_2^2\|\nbl^2v\|_2^2 +\|Q\|_2^2)
\end{split}
\end{equation*}
 Therefore,
\begin{equation*}
\begin{split}
 \int_0^T\lf(\|\Dlt\tht\|_2^2+\|\tht_z\|_{V_2}^2\rt)
 \ls& \lf|\|\tht(T)\|_{V_2}^2-\|\tht(0)\|_{V_2}^2\rt| + C\|Q\|_2^2T\\
 &+ C\int_0^T\|\nbl\tht\|_{V_2}dt +C\int_0^T\|\nbl^2v\|_2^2 dt\\
 \ls& CT^\hf +CT +CT^\hf\lf(\int_0^T\|\nbl\tht\|_{V_2}^2dt\rt)^\hf\\
 &+ C\int_0^T\|\nbl v\|_{V_1}^2dt\\
\ls& C(T^\hf+T).
\end{split}
\end{equation*}

\end{prf}

\section{Dimensions of the Global Attractor}
\label{s:dim}

 Recall the following result due to Ladyzhenskaya, see \cite{lady:90}.

\begin{thrm}
\label{t:lady}
 Let $X$ be a Hilbert space, $S : X \mapsto X$ be a map and $\A\subset X$ be
 a compact set such that $S(\A) = \A$. Suppose that there exist
 $l\in[1,+\infty)$ and $\delta\in (0, 1)$, such that $\forall a_1, a_2\in \A$,
\begin{align*}
 \|S(a_1)-S(a_2)\|_X &\ls l\|a_1-a_2\|_X ,\\
 \|Q_N \lf[S(a_1) - S(a_2)\rt]\|_X &\ls \delta\|a_1-a_2\|_X,
\end{align*}
 where $Q_N$ is the projection in $X$ onto some subspace $(X_N)^\perp$ of
 co-dimension $N\in\N$.
Then
\[ d_H(\A) \ls d_F(\A) \ls
 N\frac{\ln\lf(\frac{8Ga^2l^2}{1-\dlt^2}\rt)}{\ln\lf(\frac{2}{1+\dlt^2}\rt)},\]
 where $d_H(\A)$ and $d_F(\A)$ are the Hausdorff and fractal dimensions of
 $\A$ respectively and $Ga$ is the Gauss constant:
\[ Ga:=\frac{1}{2\pi}\beta\lf(\frac{1}{4},\frac{1}{2}\rt)
    =\frac{2}{\pi}\int_0^1\frac{dx}{\sqrt{1-x^4}}=0.8346268.... \]

\end{thrm}

 Now we use Theorem~\ref{t:lady} to prove our main result of this paper:
\begin{thrm}
\label{t:main}
 Suppose $Q\in L^2(\Omg)$ and $\alf\gs0$. Then the global attractor $\A$ has
 finite Hausdorff and fractal dimensions measured in the $V$ space.
\end{thrm}
\begin{rmrk}
 It is proved in \cite{j;t:2014} that for $Q, \prt_z Q\in L^2(\Omg)$ and
 $\alf=0$, the global attractor $\A$ has finite Hausdorff and fractal
 dimensions measured in the $V$ space. In the proof of \cite{j;t:2014},
 $H^2$ uniform boundedness on the global attractor is used. In the following,
 we give a different proof which does not use $H^2$ uniform boundedness. Thus,
 it extends the result of finite dimensionality of $\A$ to the more general
 cases.
\end{rmrk}

\begin{prf}

 Suppose $\lf(v^{(i)}, \tht^{(i)}\rt)$, $i=1,2$, are two strong
 solutions to the 3D viscous PEs, with the initial data
 $\lf(v^{(i)}_0,\tht^{(i)}_0\rt)\in\A $, where $\A$ is the global attractor
 obtained in Theorem~\ref{t:attractor}. Without losing generality, we can
 assume that $\lf(v^{(i)}_0, \tht^{(i)}_0\rt)\in D(A_1)\times D(A_2)$ for
 $i=1,2$. Therefore $\|v^{i}_t(0)\|_2, \|\tht^{i}_t(0)\|_2$ are finite for
 $i=1,2$.
 Since $\lf(v^{(i)}_0,\tht^{(i)}_0\rt)\in\A $, by Theorem~\ref{t:t}, without
 losing generality, we can assume that $\|v^{i}_t(t)\|_2, \|\tht^{i}_t(t)\|_2$
 are uniformed bounded for $t\in[0,+\infty)$, with the bounds uniform on the
 global attractor $\A$.
 Moreover, $\lf(v^{(i)}(t), \tht^{(i)}(t)\rt)\in D(A_1)\times D(A_2)$ for any
 $t\gs0$.

 Denote:
\[ u = v^{(1)}-v^{(2)}, \quad \eta = \tht^{(1)}-\tht^{(2)},
  \quad q(x,y,t) = p^{(1)}(x,y,t)-p^{(2)}(x,y,t). \]
 Notice that here $u$ is different from that in Section~\ref{s:t}.

 We have the following equations for $u$ and $\eta$:
\begin{equation*}
\begin{split}
 \Pt u  &+ L_1 u + (u\cdot\nbl)v^{(1)} + (v^{(2)}\cdot\nbl)u\\
 &-\lf(\int_{-h}^z\nbl\cdot u(x,y,\xi, t)d\xi \rt)\Pz v^{(1)}
  -\lf(\int_{-h}^z\nbl\cdot v^{(2)}(x,y,\xi,t)d\xi \rt)\Pz u\\
 &+\nbl q(x,y,t) - \int_{-h}^z\nbl\eta(x,y,\xi,t)d\xi + fu^\bot =0,
\end{split}
\end{equation*}
\begin{equation*}
\begin{split}
 \Pt\eta &+ L_2\eta + u\cdot\nbl\tht^{(1)} + v^{(2)}\cdot\nbl\eta\\
 &-\lf(\int_{-h}^z\nbl\cdot u(x,y,\xi,t)d\xi\rt)\Pz\tht^{(1)}
 -\lf(\int_{-h}^z\nbl\cdot v^{(2)}(x,y,\xi,t)d\xi\rt)\Pz\eta=0.
\end{split}
\end{equation*}

 Let $\{\lmd_{k}\}_1^\infty$ be the set of eigenvalues of $A_1$ and
 $\{\sgm_{k}\}_1^\infty$ be the set of eigenvalues of $A_2$. Let
 $P_{1,n}$ be the orthogonal projector in $H_1$ on the subspace of spanned by
 the first $n$ eigenvectors associated with $\lmd_1, \dots, \lmd_n$, and
 $P_{2,n}$ be the orthogonal projector in $H_2$ on the subspace of $H_2$
 spanned by the first $n$ eigenvectors associated with $\sgm_1, \dots, \sgm_n$.
 Let $Q_{i,n}= I - P_{i,n}$, $i=1,2$. By the classic spectral theory of compact
 operators, we know that
\[ \lim_{n\goto\infty}\lmd_n = \lim_{n\goto\infty}\sgm_n =+\infty. \]

 Now, we consider the following equations for $Q_{1,n}u$ and $Q_{2,n}\eta$:
\begin{equation}
\label{e:u.V}
\begin{split}
 \half\dt{\|Q_{1,n} u\|_{V_1}^2} +& \|A_1Q_{1,n}u\|_2^2\\
=&-\lra{(u\cdot\nbl)v^{(1)}}{A_1Q_{1,n}u}-\lra{(v^{(2)}\cdot\nbl)u}{A_1Q_{1,n}u}\\
 &+\lra{\lf(\int_{-h}^z\nbl\cdot u(x,y,\xi,t)d\xi\rt)\Pz v^{(1)}}{A_1Q_{1,n}u}\\
 &+\lra{\lf(\int_{-h}^z\nbl\cdot v^{(2)}(x,y,\xi,t)d\xi\rt)\Pz u}{A_1Q_{1,n}u}\\
 &+\lra{\int_{-h}^z\eta(x,y,\xi,t)d\xi}{A_1 Q_{1,n}u} - \lra{fu^\bot}{A_1Q_{1,n}u}\\
 &\equiv I_1 + I_2 + I_3 + I_4 + I_5 + I_6.
\end{split} 
\end{equation}
\begin{equation}
\label{e:eta.V}
\begin{split}
 \hf\dt{\|Q_{2,n}\eta\|_{V_2}^2} +& \|A_2Q_{2,n}\eta\|_2^2\\
 =& -\lra{(u\cdot\nbl)\tht^{(1)}}{A_2Q_{2,n}\eta}
  - \lra{(v^{(2)}\cdot\nbl)\eta}{A_2Q_{2,n}\eta}\\
  &+\lra{\lf(\int_{-h}^z\nbl\cdot u(x,y,\xi,t)d\xi\rt)\Pz\tht^{(1)}}{A_2Q_{2,n}\eta}\\
  &+\lra{\lf(\int_{-h}^z\nbl\cdot v^{(2)}(x,y,\xi,t)d\xi\rt)\Pz\eta}{A_2Q_{2,n}\eta}\\
  & \equiv J_1 + J_2 + J_3 + J_4
\end{split} 
\end{equation}
 First, we estimate the right-hand side of (\ref{e:u.V}) term by term.
\begin{equation*}
\begin{split}
 I_1 &\ls \|u\|_6\|\nbl v^{(1)}\|_3 \|A_1Q_{1,n}u\|_2\\
 &\ls C\|u\|_{V_1}\|v^{(1)}\|_{V_1}^\hf\|A_1 v^{(1)}\|_2^\hf \|A_1Q_{1,n}u\|_2\\
 &\ls C\|u\|_{V_1}^2 \|v^{(1)}\|_{V_1}\|A_1 v^{(1)}\|_2 
 + \frac{1}{12}\|A_1Q_{1,n} u\|_2^2.
\end{split}
\end{equation*}
 Using the Agmon inequality, we obtain
\begin{equation*}
\begin{split}
 I_2 &\ls \|v^{(2)}\|_\infty\|\nbl u\|_2\|A_1Q_{1,n}u\|_2\\
 &\ls C\|v^{(2)}\|_{V_1}^\hf\|A_1 v^{(2)}\|_2^\hf \|u\|_{V_1}\|A_1Q_{1,n}u\|_2\\
 &\ls C\|v^{(2)}\|_{V_1}\|A_1 v^{(2)}\|_2\|u\|_{V_1}^2
 + \frac{1}{12}\|A_1Q_{1,n}u\|_2^2.
\end{split}
\end{equation*}
 By Lemma~\ref{l:ju}, we have
\begin{equation*}
\begin{split}
 I_3 
\ls& C\|\nbl u\|_2^\hf\|\nbl u\|_{H^1}^\hf
      \|\prt_zv\|_2^\hf\|\prt_z v\|_{H^1}^\hf \|A_1Q_{1,n}u\|_2\\
\ls& C \|v_z^{(1)}\|_2\|v_z^{(1)}\|_{H^1} \|u\|_{V_1} \|A_1u\|_2
   + \frac{1}{12} \|A_1Q_{1,n}u\|_2^2.\\
\end{split}
\end{equation*}
 Similarly, we have
\begin{equation*}
\begin{split}
 I_4
\ls& C\|\nbl v^{(2)}\|_2^\hf \|\nbl^2 v^{(2)}\|_2^\hf
   \|u_z\|_2^\hf \|\nbl u_z\|_2^\hf\|A_1Q_{1,n}u\|_2\\
\ls& C\|\nbl v^{(2)}\|_2\|\nbl^2 v^{(2)}\|_2\|u_z\|_2\|\nbl u_z\|_2
   + \frac{1}{12} \|A_1Q_{1,n}u\|_2^2.
\end{split}
\end{equation*}
 The last two terms can be estimated as follows:
\begin{equation*}
 I_5 \ls C \|\eta\|_{V_2} \|A_1Q_{1,n}u\|_2
\ls C\|\eta\|_{V_2}^2 + \frac{1}{12}\|A_1Q_{1,n}u\|_2^2,
\end{equation*}
\begin{equation*}
 I_6 \ls C \|u\|_{V_1} \|A_1Q_{1,n}u\|_2
 \ls C\|u\|_{V_1}^2 + \frac{1}{12}\|A_1Q_{1,n}u\|_2^2.
\end{equation*}
 Next, we estimate the right-hand side of (\ref{e:eta.V}) term by term.
 Similar to the estimate for $I_1$, we have
\begin{equation*}
\begin{split}
 J_1 &\ls \|u\|_6\|\nbl \tht^{(1)}\|_3 \|A_2Q_{2,n}\eta\|_2\\
 &\ls C\|u\|_{V_1}\|\tht^{(1)}\|_{V_2}^\hf\|A_2 \tht^{(1)}\|_2^\hf
 \|A_2Q_{2,n}\eta\|_2\\
 &\ls C\|\tht^{(1)}\|_{V_2}\|A_2\tht^{(1)}\|_2
 \|u\|_{V_1}^2 + \frac{1}{8}\|A_2Q_{2,n}\eta\|_2^2.
\end{split}
\end{equation*}
 Similar to the estimate for $I_2$, we have
\begin{equation*}
\begin{split}
 J_2
&\ls \|v^{(2)}\|_\infty\|\nbl\eta\|_{V_2} \|A_2Q_{2,n}\eta\|_2\\
&\ls C\|v^{(2)}\|_{V_1}^\hf\|A_1 v^{(2)}\|_2^\hf\|\eta\|_{V_2}
     \|A_2Q_{2,n}\eta\|_2\\
&\ls C\|v^{(2)}\|_{V_1}\|A_1 v^{(2)}\|_2\|\eta\|_{V_2}^2
 + \frac{1}{8}\|A_2Q_{2,n}\eta\|_2^2.
\end{split}
\end{equation*}
 Similar to the estimates for $I_3$ and $I_4$, we have
\begin{equation*}
\begin{split}
 J_3
\ls& C\|\nbl u\|_2^\hf\|\nbl^2u\|_2^\hf
      \|\tht^{(1)}_z\|_2^\hf\|\nbl\tht^{(1)}_z\|_2^\hf \|A_2Q_{2,n}\eta\|_2\\
\ls& C\|\tht^{(1)}_z\|_2 \|\nbl\tht^{(1)}_z\|_2 \|\nbl u\|_2 \|\nbl^2u\|_2
 + \frac{1}{8}\|A_2Q_{2,n}\eta\|_2^2,
\end{split}
\end{equation*}
 and
\begin{equation*}
\begin{split}
 J_4
\ls& C\|v^{(2)}\|_{V_1}^\hf \|A_1 v^{(1)}\|_2^\hf
 \|\eta_z\|_2^\hf \|\eta_z\|_{H^1}^\hf\|A_2Q_{2,n}\eta\|_2\\
\ls& C\|v^{(2)}\|_{V_1}\|A_1 v^{(2)}\|_2\|\eta_z\|_2\|\eta_z\|_{H^1}
 +\frac{1}{8} \|A_2Q_{2,n}\eta\|_2^2\\
\ls& C\|v^{(2)}\|_{V_1}\|A_1 v^{(2)}\|_2\|\eta\|_{V_2}\|A_2\eta\|_2
 +\frac{1}{8} \|A_2Q_{2,n}\eta\|_2^2.
\end{split}
\end{equation*}
 Now adding \eqref{e:u.V} and \eqref{e:eta.V} and using the estimates for
 $I_i$ ($i=1,...,6$) and $J_j$ ($j=1,...,4$), and denoting
\begin{align*}
 y(t):= & \|u(t)\|_{V_1}^2 + \|\eta(t)\|_{V_2}^2,\\
 z(t):= & \|Q_{1,n}u(t)\|_{V_1}^2 + \|Q_{2,n}\eta(t)\|_{V_2}^2,
\end{align*}
 we have, noticing that $(v^{(i)}(t),\tht^{(i)}(t))\in\A$, for $i=1,2$,
\begin{equation}
\label{e:z.0}
\begin{split}
 z'(t) +& \|A_1Q_{1,n}u\|_2^2 + \|A_2Q_{2,n}\eta\|_2^2\\
 \ls& Cy(t)(1+\|A_1v^{(1)}\|_2+\|A_1v^{(2)}\|_2+\|A_2\tht^{(1)}\|_2)\\
 &+ C\|u\|_{V_1}\|A_1 u\|_2
   (\|\nbl v_z^{(1)}\|_2+\|\nbl\tht^{(1)}_z\|_2 +\|\nbl^2v^{(2)}\|_2)\\
 &+ C\|\eta\|_{V_1}\|A_2\eta\|_2\|\nbl^2v^{(2)}\|_2,
\end{split}
\end{equation}
 where $C$ is a generic positive constant independent of $t$, $n$ and the
 initial data. Denote $\rho_n:=\min(\lmd_n^2,\sgm_n^2)$. Then, it follows
 from \eqref{e:z.0} that
\begin{equation}
\label{e:z.1}
\begin{split}
 z'(t)+\rho_nz(t)
\ls & Cy(t)(1+\|A_1v^{(1)}\|_2+\|A_1v^{(2)}\|_2+\|A_2\tht^{(1)}\|_2)\\
   &+ C y^\hf(t)(\|A_1u\|_2^2+\|A_2\eta\|_2^2)^\hf\\
    &\times(\|\nbl v_z^{(1)}\|_2+\|\nbl\tht^{(1)}_z\|_2 +\|\nbl^2v^{(2)}\|_2).
\end{split}
\end{equation}
 Now, integrating \eqref{e:z.1} with respect to $t\in[0,T]$, we obtain
\begin{equation*}
\begin{split}
 z(T) \ls & e^{-\rho_n T}z(0) + Ce^{-\rho_n T}\int_0^T e^{\rho_n t}y(t)dt\\
          & + Ce^{-\rho_n T}\int_0^T e^{\rho_n t}y(t)
 (\|A_1v^{(1)}\|_2+\|A_1v^{(2)}\|_2+\|A_2\tht^{(1)}\|_2)dt\\
 & +Ce^{-\rho_n T}\int_0^T e^{\rho_n t}y^\hf(t)
       (\|A_1 u(t)\|_2^2+\|A_2\eta(t)\|_2^2)^\hf\\
 &\times(\|\nbl v_z^{(1)}\|_2+\|\nbl\tht^{(1)}_z\|_2 +\|\nbl^2v^{(2)}\|_2) dt\\
   =:& Z_1 +Z_2 +Z_3 +Z_4.
\end{split}
\end{equation*}
 First, notice that
\[ Z_1 \ls e^{-\rho_n T}y(0). \]
 Next, recall that it is proved in \cite{j:pe} that for strong solutions
 $(v^{(i)}(t),\tht^{(i)}(t))\in V$, $i=1,2$, there exists a positive continuous
 non-decreasing function $K(t)$, independent of the initial data, such that
\begin{equation}
\label{e:lpc.s}
 y(t) + \int_0^t(\|A_1 u(\tau)\|_2^2+\|A_2\eta(\tau)\|_2^2)d\tau
 \ls K(t) y(0), \quad\forall t\gs 0.
\end{equation}
 Therefore,
\[Z_2\ls CK(T)y(0)e^{-\rho_n T}\int_0^T e^{\rho_n t} dt \ls CK(T)\rho_n^{-1}y(0)\]
 and
\begin{equation*}
\begin{split}
 Z_3 \ls&  CK(T)y(0)e^{-\rho_n T}\int_0^T e^{\rho_n t}
  (\|A_1v^{(1)}\|_2+\|A_1v^{(2)}\|_2+\|A_2\tht^{(1)}\|_2)dt\\
  \ls& CK(T)y(0)e^{-\rho_n T}\lf(\int_0^T e^{2\rho_n t} dt\rt)^\hf\\
  &\times \lf[\int_0^T
   \lf(\|A_1v^{(1)}\|_2^2+\|A_1v^{(2)}\|_2^2+\|A_2\tht^{(1)}\|_2^2\rt)dt\rt]^\hf\\
  \ls& CK(T)y(0)e^{-\rho_n T}\lf(\int_0^T e^{2\rho_n t} dt\rt)^\hf
  \ls CK(T)\rho_n^{-\hf}y(0),
\end{split}
\end{equation*}
 by the uniform boundedness of
$\int_0^T
  \lf(\|A_1v^{(1)}\|_2^2+\|A_1v^{(2)}\|_2^2+\|A_2\tht^{(1)}\|_2^2\rt)dt$
 on the global attractor $\A$.

 The estimate of $Z_4$ is more complicated and needs a special treatment.
\begin{equation*}
\begin{split}
 Z_4 &\ls  CK^\hf(T)y^\hf(0)
 \int_0^T(\|A_1u(t)\|_2^2+\|A_2\eta(t)\|_2^2)^\hf\\
  &\times(\|\nbl v_z^{(1)}\|_2+\|\nbl\tht^{(1)}_z\|_2 +\|\nbl^2v^{(2)}\|_2) dt \\
 & \ls CK^\hf(T)y^\hf(0) \lf(\int_0^T\|A_1u\|_2^2+\|A_2\eta\|_2^2 dt\rt)^\hf\\
  &\times\lf[\int_0^T
  \lf(\|\nbl v_z^{(1)}\|_2^2+\|\nbl\tht^{(1)}_z\|_2^2 +\|\nbl^2v^{(2)}\|_2^2\rt)
  dt \rt]^\hf\\
 & \ls CK(T)(T^\hf+T)^\hf y(0),
\end{split}
\end{equation*}
 where we have used \eqref{e:lpc.s} and uniform continuity property as proved
 in Theorem~\ref{t:uc}.

 Combining the estimates of the $Z_i$'s, $i=1,\dots,4$, we obtain
\[ z(T) \ls \lf[e^{-\rho_n T} + C_1K(T)\rho_n^{-1}
  + C_2K(T)\rho_n^{-1/2}+C_3K(T)(T^\hf+T)^\hf\rt]y(0), \]
 where $C_i$'s, $i=1,2,3$, are all constants independent of $T$, $n$ and
 initial data.

 Notice that
\[ \lim_{T\goto 0+} K(T)(T^\hf+T)^\hf=0. \]
 Therefore, for any $\dlt\in(0,1)$, we can choose $T>0$ sufficiently small
 and uniform for any pair of $\lf(v^{(i)}, \tht^{(i)}\rt)\in\A$, $i=1,2$, such
 that
\[ C_3K(T)(T^\hf+T)^\hf \ls \frac{\dlt}{2}, \]
 and then we choose $n$ sufficiently large, which is uniform for any pair of
 $\lf(v^{(i)}, \tht^{(i)}\rt)\in\A$, $i=1,2$, such that
\[ e^{-\rho_n T} + C_1K(T)\rho_n^{-1}
  + C_2K(T)\rho_n^{-1/2} \ls \frac{\dlt}{2}. \]
 Thus, 
\[ z(T) \ls \dlt y(0). \]
 This proves the expected squeezing property of $S(T)$ for any $T>0$ that is
 sufficiently small.

 Notice that \eqref{e:lpc.s} gives the Lipschitz continuity of the solution
 map $S(t)(v_0,\tht_0):= (v(t),\tht(t))$ in the space $V$ for any fixed
 $t\in[0,\infty)$. Thus, Theorem~\ref{t:main} follows immediately from
 Theorem~\ref{t:lady}.

\end{prf}

\begin{rmrk}
 It is easy to see, by a continuation argument, that the squeezing property
 of $S(T)$ is indeed still valid for any $T>0$.
\end{rmrk}

\section{The case with ``physical boundary conditions''}
\label{s:pbc}

 The previous sections of this paper were finished in December 2014 and the
 original manuscript was submitted for publication by then. This short section
 is now added to include the case with ``physical boundary conditions''. This
 is the case that the velocity field $v$ satisfies, instead of \eqref{e:bc.v},
 the following boundary conditions:
\begin{equation}
\label{e:bc.v.p}
 (\alf_v v + v_z)|_{z=0}= v|_{(x,y)\in\prt D}= v|_{z=-h}=0,
\end{equation}
 where $\alf_v$ is a non-negative real constant independent of $v$, to be
 distinguished from the non-negative real constant $\alf$ in \eqref{e:bc.t}.

 It can be shown that {\em the main results of this paper are still valid for
 this case}. The detailed proofs are very similar to those given in the previous
 sections and thus are omitted here. However, it is {\em not} clear if the
 method of \cite{j;t:2014} can be used for this case. This is an additional
 example showing the advantage of the new approach as given in this paper when
 compared with the one presented in our previous \cite{j;t:2014}. \\

{\bf Acknowledgment:}
 This is a continuation of the recent work \cite{j;t:2014}. The author wishes
 to thank Prof. Roger Temam for his precious support of this research project.

\end{document}